\documentclass[12pt]{article}
\usepackage{amssymb,amsmath}
\usepackage{hyperref}
\usepackage{graphicx}
\usepackage{color}
\usepackage{chngcntr}

\begin{document}

\title{\LARGE\bf An iteration procedure for a two-term Machin-like formula for pi with \\
small Lehmer's measure}

\author{
\normalsize\bf S. M. Abrarov\footnote{\scriptsize{Dept. Earth and Space Science and Engineering, York University, Toronto, Canada, M3J 1P3.}}\, and B. M. Quine$^{*}$\footnote{\scriptsize{Dept. Physics and Astronomy, York University, Toronto, Canada, M3J 1P3.}}}

\date{July 26, 2017}
\maketitle

\begin{abstract}
\vspace{0.1cm}
In this paper we present a two-term Machin-like formula for pi \[\frac{\pi}{4} = 2^{k - 1}\arctan\left(\frac{1}{u_1}\right) + \arctan\left(\frac{1}{u_2}\right)\] with small Lehmer's measure $e \approx \,0.245319$ and describe iteration procedure for simplified determination of the required rational number $u_2$ at $k = 27$ and $u_1 = 85445659$. With these results we obtained a formula that has no irrational numbers involved in computation and provides $16$ digits of pi at each increment by one of the summation terms. This is the smallest Lehmer's measure ever reported for the Machin-like formulas for pi.
\vspace{0.25cm}
\\
\noindent {\bf Keywords:} Machin-like formula, constant pi, Lehmer's measure, arctangent function
\vspace{0.25cm}
\end{abstract}

\section{Introduction}

Machin-like formulas for pi can be generalized in form \cite{Lehmer1938, Borwein1987}
\begin{equation}\label{eq_1}
\frac{\pi }{4}=\sum\limits_{k=1}^{K}{{{\alpha }_{k}}}\arctan \left( \frac{1}{{{\beta }_{k}}} \right),
\end{equation}
where ${{\alpha }_{k}}$ and ${{\beta }_{k}}$ are integers or rationals. Since the Maclaurin series expansion of the arctangent function is given by
$$
\arctan \left( x \right)=x-\frac{{{x}^{3}}}{3}+\frac{{{x}^{5}}}{5}-\frac{{{x}^{7}}}{7}+\ldots \,\,,
$$
we can rewrite it in a simplified form as
$$
\arctan \left( x \right)=x+O\left( {{x}^{3}} \right).
$$
Consequently, due to vanishing term $O\left( {{x}^{3}} \right)$ one can expect a rapid improvement in accuracy of the arctangent function as its argument $x$ decreases by absolute value. Therefore, it is more efficient to compute pi when arguments $1/{{\beta }_{k}}$ of the arctangent function in equation \eqref{eq_1} are smaller by absolute value.

The measure defined as
\begin{equation}\label{eq_2}
e=\sum\limits_{k=1}^{K}{\frac{1}{{{\log }_{10}}\left( \left| {{\beta }_{k}} \right| \right)}}
\end{equation}
can be used to quantify the computational efficiency of the Machin-like formulas \eqref{eq_1} for pi. In particular, Lehmer in his paper \cite{Lehmer1938} stated that the measure \eqref{eq_2} shows how much labor is required for a specific Machin-like formula to compute pi. From this statement it follows that the computational efficiency of the given Machin-like formula for pi is higher when its measure \eqref{eq_2} is smaller. As we can see, the Lehmer's measure \eqref{eq_2} decreases at smaller number of the arctangent function terms $K$ and at larger values of the parameters $\beta_k$ by absolute value in the equation \eqref{eq_1}.

In $2002$ Kanada while breaking a record applied the following self-checking pair of the Machin-like formulas
\small
\[
\frac{\pi }{4}=44\arctan \left( \frac{1}{57} \right)+7\arctan \left( \frac{1}{239} \right)-12\arctan \left( \frac{1}{682} \right)+24\arctan \left( \frac{1}{12943} \right) 
\]
\normalsize
and
\small
\[
\frac{\pi }{4}=12\arctan \left( \frac{1}{49} \right)+32\arctan \left( \frac{1}{57} \right)-5\arctan \left( \frac{1}{239} \right)+12\arctan \left( \frac{1}{110443} \right) 
\]
\normalsize
that enabled him to compute the constant pi with correct number of digits exceeding one trillion \cite{Calcut2009}. This signifies a strong potential of the Machin-like formulas \eqref{eq_1} even at relatively large Lehmer's measure. Therefore, the derivation of the Machin-like formulas for pi with reduced Lehmer's measure remains an interesting topic and many new identities have been reported in the modern literature \cite{Wetherfield1996, Wetherfield1997, Chien-Lih1997, Nimbran2010, Arndt2011, Arndt} by using, for example, the Todd's process described in the work \cite{Todd1949}.

Since the identity
\[
\pi/4=\arctan\left(1\right)
\]
is the only Machin-like formula for pi with a single term \cite{Calcut2009}, this fact motivated us to develop a new methodology to gain computational efficiency by minimizing in equation \eqref{eq_1} the number of terms to two with smaller arguments (by absolute value) of the arctangent function \cite{Abrarov2017}.

In this paper we describe iteration procedure that can be used to simplify significantly the computation of a two-term Machin-like formula for pi with small Lehmer's measure. This approach leads to a rapidly convergent formula for pi consisting of the rational numbers only. In particular, the computational test demonstrates that each increment of the summation terms just by one contributes to $16$ additional digits of the constant pi. We also show that the iteration method in determination of the expansion coefficients excludes all complex numbers in computation. The absence of the irrational and complex numbers as well as the rapid convergence and simplicity of the proposed formula may be promising in computing pi. Due to no any theoretical restrictions, the convergence rate of this formula for pi can be increased further. This is practically feasible since the Lehmer's measure decreases with increasing the integer $k$.

\section{Algorithmic implementation}

\subsection{Iteration procedure}

In our previous publication we have shown that the following equation \cite{Abrarov2017}
$$
\frac{\pi }{4}={{2}^{k-1}}\arctan \left( \frac{\sqrt{2-{{a}_{k-1}}}}{{{a}_{k}}} \right),
$$
where 
$$
{{a}_{k}}=\sqrt{2+{{a}_{k-1}}}, \quad {{a}_{1}}=\sqrt{2},
$$
can be rewritten as a two-term Machin-like formula for pi as given by
\begin{equation}\label{eq_3}
\frac{\pi }{4}={{2}^{k-1}}\arctan \left( \frac{1}{{{u}_{1}}} \right)+\arctan \left( \frac{1}{{{u}_{2}}} \right),
\end{equation}
where ${{u}_{1}}$ is a positive rational number such that
\begin{equation}\label{eq_4}
{{u}_{1}}=\frac{{{a}_{k}}}{\sqrt{2-{{a}_{k-1}}}}+\varepsilon ,	\qquad{{u}_{1}}>>\left| \varepsilon  \right|
\end{equation}
and
\begin{equation}\label{eq_5}
{{u}_{2}}=\frac{2}{{{\left( \left( {{u}_{1}}+i \right)/\left( {{u}_{1}}-i \right) \right)}^{{{2}^{k-1}}}}-i}-i.
\end{equation}
Since
$$
\frac{1}{{{u}_{2}}}=\frac{2}{{{\left( \left( {{u}_{1}}+i \right)/\left( {{u}_{1}}-i \right) \right)}^{{{2}^{k-1}}}}+i}+i
$$
the equation \eqref{eq_3} can also be represented in form
\begin{equation}\label{eq_6}
\frac{\pi }{4}={{2}^{k-1}}\arctan \left( \frac{1}{{{u}_{1}}} \right)+\arctan \left( \frac{2}{{{\left( \left( {{u}_{1}}+i \right)/\left( {{u}_{1}}-i \right) \right)}^{{{2}^{k-1}}}}+i}+i \right).
\end{equation}

With equations \eqref{eq_3}, \eqref{eq_4} and \eqref{eq_5} it is very easy to derive the well-known formula for pi that was originally discovered in $1706$ by English mathematician John Machin and named in his honor \cite{Arndt2001}. Particularly, at $k=3$ we have
$$
\frac{{{a}_{3}}}{\sqrt{2-{{a}_{2}}}}=\frac{\sqrt{2+\sqrt{2+\sqrt{2}}}}{\sqrt{2-\sqrt{2+\sqrt{2}}}}=\text{5}\text{.02733949212584810451}\ldots \,\,\left( \text{irrational} \right).
$$
Assuming that the error term
\[
\varepsilon =-0.\text{02733949212584810451}\ldots \,\,\left( \text{irrational} \right)
\]
from equation \eqref{eq_4} it follows that ${{u}_{1}}=5$. Consequently, substituting $k=3$ and ${{u}_{1}}=5$ into equation \eqref{eq_5} we can readily find that
$$
{u_2}=\frac{2}{{{\left( \left( {5}+i \right)/\left( {5}-i \right) \right)}^{{{2}^{3-1}}}}-i}-i=-{239}.
$$
Using $k=3$ and $u_2=-239$ in equation \eqref{eq_3} immediately yields the original Machin's formula for pi \cite{Lehmer1938, Borwein1987, Arndt2001}
$$
\begin{aligned}
\frac{\pi }{4}&=4\arctan \left( \frac{1}{5} \right)+\arctan \left( \frac{1}{-239} \right) \\ 
& =4\arctan \left( \frac{1}{5} \right)-\arctan \left( \frac{1}{239} \right).  
\end{aligned}
$$

Although computation of the rational number ${{u}_{2}}$ is simple, the application of the formula \eqref{eq_5} leads to some complexities due to rapidly growing power ${{2}^{k-1}}$ as the integer $k$ increases. For example, at $k=27$ the value of the power becomes colossal ${{2}^{k-1}}=67108864$. As a result, the determination of the value ${{u}_{2}}$ by straightforward application of equation \eqref{eq_5} requires extended computer memory usage and becomes extremely time-consuming. However, these computational complexities can be successively resolved by applying the iteration procedure that we developed in this work.

Defining the real and imaginary parts as
\[
\tag{7a}\label{eq_7a}
x_1=\operatorname{Re}\left[ \frac{{{u}_{1}}+i}{{{u}_{1}}-i} \right]=\frac{u_{1}^{2}-1}{u_{1}^{2}+1}
\]
and
\[
\tag{7b}\label{eq_7b}
y_1=\operatorname{Im}\left[ \frac{{{u}_{1}}+i}{{{u}_{1}}-i} \right]=\frac{2u_1}{u_{1}^{2}+1},
\]
respectively, the equation \eqref{eq_5} can be conveniently rewritten in form
\setcounter{equation}{7}
\begin{equation}\label{eq_8}
{{u}_{2}}=\frac{2}{{{\left( x_1+iy_1 \right)}^{{{2}^{k-1}}}}-i}-i.
\end{equation}

It is not difficult to see by induction that
\small
\begin{equation}\label{eq_9}
\begin{aligned}
\left( {{x}_{1}}\right. &+ \left. i{{y}_{1}} \right)^{{{2}^{k-1}}}\\
&=\overbrace{{{\left( {{\left( {{\left( {{\left( {{x}_{1}}+i{{y}_{1}} \right)}^{2}} \right)}^{2}} \right)}^{2\,\,\,\cdots }} \right)}^{2}}}^{k-1\,\,\text{powers}\,\,\text{of}\,\text{2}}=\overbrace{{{\left( {{\left( {{\left( {{\left( {{x}_{2}}+i{{y}_{2}} \right)}^{2}} \right)}^{2}} \right)}^{2\,\,\,\cdots }} \right)}^{2}}}^{k-2\,\,\text{powers}\,\,\text{of}\,\text{2}} \\ 
 & =\overbrace{{{\left( {{\left( {{\left( {{\left( {{x}_{3}}+i{{y}_{3}} \right)}^{2}} \right)}^{2}} \right)}^{2\,\,\,\cdots }} \right)}^{2}}}^{k-3\,\,\text{powers}\,\,\text{of}\,\text{2}}=\cdots =\overbrace{{{\left( {{\left( {{\left( {{\left( {{x}_{n}}+i{{y}_{n}} \right)}^{2}} \right)}^{2}} \right)}^{2\,\,\,\cdots }} \right)}^{2}}}^{k-n\,\,\text{powers}\,\,\text{of}\,\text{2}}=\cdots  \\ 
 & ={{\left( {{\left( {{x}_{k-2}}+i{{y}_{k-2}} \right)}^{2}} \right)}^{2}}={{\left( {{x}_{k-1}}+i{{y}_{k-1}} \right)}^{2}}={{x}_{k}}+i{{y}_{k}},  
\end{aligned}
\end{equation}
\normalsize
where the numbers ${{x}_{n}}$ and ${{y}_{n}}$ can be found by the following iteration procedure
\begin{equation}\label{eq_10}
\left\{ {
\begin{aligned}
&x_n=x_{n-1}^2-y_{n-1}^2 \\
&y_n=2x_{n-1}y_{n-1}, \qquad n=\left\{2,3,4,\,\ldots\,,k\right\}.
\end{aligned}
}\right.
\end{equation}
Consequently, from the equations \eqref{eq_8} and \eqref{eq_9} it follows that
\begin{equation}\label{eq_11}
{{u}_{2}}=\frac{2}{{{x}_{k}}+i{{y}_{k}}-i}-i=\frac{2{{x}_{k}}}{x_{k}^{2}+{{\left( {{y}_{k}}-1 \right)}^{2}}}+i\left( \frac{2\left( 1-{{y}_{k}} \right)}{x_{k}^{2}+{{\left( {{y}_{k}}-1 \right)}^{2}}}-1 \right).
\end{equation}

\subsubsection*{Theorem 1}

The value ${{u}_{2}}$ is real if ${{u}_{1}}$ is real.

\subsubsection*{Proof}

Using de Moivre's formula we can write the complex number in polar form as follows
\footnotesize
$$
{{\left( x_1+iy_1 \right)}^{{{2}^{k-1}}}}={{\left( {x_1^2}+{y_1^2} \right)}^{{{2}^{k-2}}}}\left( \cos \left( {{2}^{k-1}}\text{Arg}\left( x_1+iy_1 \right) \right)+i\sin \left( {{2}^{k-1}}\text{Arg}\left( x_1+iy_1 \right) \right) \right)
$$
\normalsize
Consequently, applying this identity into equation \eqref{eq_8} and representing $x_1$ and $y_1$ in accordance with equations \eqref{eq_7a} and \eqref{eq_7b}, respectively, after some trivial rearrangements we obtain
\begin{equation}\label{eq_12}
{{u}_{2}}=\frac{\cos \left( {{2}^{k-1}}\text{Arg}\left( \frac{{{u}_{1}}+i}{{{u}_{1}}-i} \right) \right)}{1-\sin \left( {{2}^{k-1}}\text{Arg}\left( \frac{{{u}_{1}}+i}{{{u}_{1}}-i} \right) \right)}.
\end{equation}
From equation \eqref{eq_4} it follows that
$$
{{u}_{1}}\ge \frac{\sqrt{2+\sqrt{2}}}{\sqrt{2-\sqrt{2}}}+\varepsilon \Rightarrow {{u}_{1}}>1.
$$ 
As a consequence, the following inequality
$$
\operatorname{Re}\left[{\frac{{{u}_{1}}+i}{{{u}_{1}}-i}}\right]=\frac{{{u}_{1}^2}-1}{u_{1}^{2}+1}>0
$$
is satisfied to validate the relation between the principal value argument and the arctangent function as given by
\[
\begin{aligned}
\text{Arg}\left( \frac{{{u}_{1}}+i}{{{u}_{1}}-i} \right)&=\text{Arg}\left( \frac{u_{1}^{2}-1}{u_{1}^{2}+1}+i\frac{2{{u}_{1}}}{u_{1}^{2}+1} \right) \\ 
& =\arctan \left( \left( \frac{2{{u}_{1}}}{u_{1}^{2}+1} \right)/\left( \frac{u_{1}^{2}-1}{u_{1}^{2}+1} \right) \right)=\arctan \left( \frac{2{{u}_{1}}}{u_{1}^{2}-1} \right).
\end{aligned}
\]
Applying this relation, the equation \eqref{eq_12} can be simplified and represented as
\begin{equation}\label{eq_13}
{{u}_{2}}=\frac{\cos \left( {{2}^{k-1}}\arctan \left( \frac{2{{u}_{1}}}{u_{1}^{2}-1} \right) \right)}{1-\sin \left( {{2}^{k-1}}\arctan \left( \frac{2{{u}_{1}}}{u_{1}^{2}-1} \right) \right)}.
\end{equation}
As we can see now the value ${{u}_{2}}$ is purely real since $u_1$ is real. This signifies that the imaginary part in the equation \eqref{eq_11} must be zero. This completes the proof.

\subsubsection*{Corollary}

Since the imaginary part of the equation \eqref{eq_11} is zero, it follows that
$$
\frac{2\left( 1-{{y}_{k}} \right)}{x_{k}^{2}+{{\left( {{y}_{k}}-1 \right)}^{2}}}=1\Leftrightarrow x_{k}^{2}=1-{{y}_{k}^2}.
$$
Consequently, the equation (11) can be greatly simplified as
\begin{equation}\label{eq_14}
{{u}_{2}}=\frac{{{x}_{k}}}{1-{{y}_{k}}}.
\end{equation}

\vspace{0.25cm}
In our previous publication we have shown already that the value ${{u}_{2}}$ must be rational when the value ${{u}_{1}}$ is rational \cite{Abrarov2017}. Here we show an alternative proof based on iteration.

\subsubsection*{Theorem 2}

The value ${{u}_{2}}$ is rational if ${{u}_{1}}$ is rational.

\subsubsection*{Proof}

Since the value ${{u}_{1}}$ is rational, the values $x_1$ and $y_1$ are also rationals as it follows from the equations \eqref{eq_7a} and \eqref{eq_7b}. This signifies that all intermediate values ${{x}_{n}}$ and ${{y}_{n}}$ obtained by iteration with help of the set \eqref{eq_10} are also rationals. Therefore, from equations \eqref{eq_9}, \eqref{eq_10} and \eqref{eq_14} it follows that ${{u}_{2}}$ must be rational since ${{x}_{k}}$ and ${{y}_{k}}$ are both rationals. This completes the proof.

\subsection{Numerical results}

At each successive step of iteration the number of the digits in ${{x}_{n}}$ and ${{y}_{n}}$ considerably increases. Therefore, with a typical desktop computer we could perform previously the computations up to $k=23$ only \cite{Abrarov2017}. In this work we applied a supercomputer provided at the Algonquin Radio Observatory, Canada. This enabled us to increase significantly the integer $k$ up to $27$.

According to the iteration procedure discussed above the rational number ${{u}_{2}}$ can be computed by using a simplified variation of the formula \eqref{eq_14} as follows
$$
{{u}_{2}}=\frac{\text{num}\left( {{x}_{k}} \right)}{\text{den}\left( {{y}_{k}} \right)-\text{num}\left( {{y}_{k}} \right)},
$$
where the notations $\text{num}\left( \ldots  \right)$ and $\text{den}\left( \ldots  \right)$ denote the numerator and denominator, respectively. This simplification is possible since two values ${{x}_{k}}$ and ${{y}_{k}}$ have same denominator.

At $k=27$ we obtain
$$
\frac{{{a}_{27}}}{\sqrt{2-{{a}_{26}}}}=8.54456594470539448216\ldots \times {{10}^{7}}\,\,\left( \text{irrational} \right).
$$
We can choose the error term to be
$$
\varepsilon =-0.00000004470539448216\ldots \times {{10}^{7}}\,\,\left( \text{irrational} \right).
$$
Consequently, the rational number is ${{u}_{1}}=85445659$.

Lastly, using the described iteration procedure we obtain
\begin{equation}\label{eq_15}
\begin{aligned}
{{u}_{2}}&=-\frac{\overbrace{\text{2368557598}\ldots \text{9903554561}}^{\text{522,185,816}\,\,\text{digits}}}{\underbrace{\text{9732933578}\ldots \text{4975692799}}_{\text{522,185,807}\,\,\text{digits}}} \\ 
 & =-\text{2}\text{.43354953523904089818}\ldots \times {{10}^{8}}\,\,\left( \text{rational} \right).  
\end{aligned}
\end{equation}

The interested reader can download the computed rational number $u_2$ with all digits in the numerator and denominator \cite{u2.txt}.

\subsection{Lehmer's measure}

In 1938 Lehmer in his paper \cite{Lehmer1938} showed the three-term Machin-like formula for pi\,\footnote{In fact, Lehmer also suggested to reduce the measure $e$ from each term proportional to $\arctan\left(\frac{1}{10^q}\right)$, where $q$ is a positive integer. However, Lehmer implied these reductions only for manual calculations and ruled them out if computer is applied. Therefore, without these reductions the value $e \approx 1.5279$ is the smallest measure shown in the paper \cite{Lehmer1938}.}
$$
\frac{\pi}{4} =22\arctan\left(\frac{1}{26}\right)-2\arctan\left(\frac{1}{2057}\right)-5\arctan\left(\frac{38479}{3240647}\right),
$$
with perhaps the smallest measure $e \approx 1.5279$ known by that time. However, applying equations \eqref{eq_3}, \eqref{eq_4} and \eqref{eq_5} it is not difficult to derive the two-term Machin-like formula for pi with Lehmer's measure $e$ less than this value even at relatively small integer $k$. Let's take, for example, $k=6$ and since
$$
\begin{aligned}
\frac{a_6}{\sqrt{2-a_5}}&=\frac{\sqrt{2+\sqrt{2+\sqrt{2+\sqrt{2+\sqrt{2+\sqrt{2}}}}}}}{\sqrt{2-\sqrt{2+\sqrt{2+\sqrt{2+\sqrt{2+\sqrt{2}}}}}}} \\
&=40.73548387208330180074\ldots \, \left(\text{irrational}\right)
\end{aligned}
$$
we can choose the error term to be
$$
\epsilon = -0.73548387208330180074\ldots \, \left(\text{irrational}\right).
$$
Consequently, from equation \eqref{eq_4} we can find that $u_1=40$. Substituting now $k=6$ and $u_1=40$ into equations \eqref{eq_5} and then \eqref{eq_3} we get the two-term Machin-like formula for pi
\footnotesize
\[
\begin{aligned}
\frac{\pi }{4}=&\,{{2}^{6-1}}\arctan \left( \frac{1}{40} \right)\\
&+\arctan \left( -\frac{38035138859000075702655846657186322249216830232319}{2634699316100146880926635665506082395762836079845121} \right)
\end{aligned}
\]
\normalsize
or
\footnotesize
\[
\begin{aligned}
\frac{\pi }{4}=&\,32\arctan \left( \frac{1}{40} \right)\\
&-\arctan \left(\frac{38035138859000075702655846657186322249216830232319}{2634699316100146880926635665506082395762836079845121} \right)
\end{aligned}
\]
\normalsize
with Lehmer's measure $e \approx 1.16751$ only. The following Mathematica code:
\small
\begin{verbatim}
32*ArcTan[1/40] - 
  ArcTan[38035138859000075702655846657186322249216830232319/
    2634699316100146880926635665506082395762836079845121] == Pi/4
\end{verbatim}
\normalsize
returns the output {\ttfamily{True}}\,\footnote{This code verifies whether or not the left side of the equation is equal to its right side.}.

At $k=27$, ${{u}_{1}}=85445659$ and corresponding ${{u}_{2}}$ (see equation \eqref{eq_15}) the Lehmer's measure \eqref{eq_2} for the two-term Machin-like formula \eqref{eq_3} for pi becomes $e\approx \,0.245319$. To the best of our knowledge this value of the Lehmer's measure is the smallest ever reported in scientific literature for the Machin-like formulas \eqref{eq_1} for pi. For example, one of the smallest known nowadays Lehmer's measure $e\approx 1.51244$ corresponds to the following six-term Machin-like formula for pi \cite{Chien-Lih1997}
\footnotesize
\[
\begin{aligned}
\frac{\pi }{4}=& 183\arctan \left( \frac{1}{239} \right)+32\arctan \left( \frac{1}{1023} \right)-68\arctan \left( \frac{1}{5832} \right)\\
&+12\arctan \left( \frac{1}{110443} \right)-12\arctan \left( \frac{1}{4841182} \right)-100\arctan \left( \frac{1}{6826318} \right).  
\end{aligned}
\]
\normalsize
As we can see, the obtained Lehmer's measure is about $6$ times smaller than that of corresponding to the six-term Machin-like formula for pi. Since the smaller Lehmer's measure characterizes the higher computational efficiency, the proposed two-term Machin-like formula for pi may be promising for computation of the constant pi.

\section{Convergence}

In our earlier publication we have derived a new formula for the arctangent function \cite{Abrarov2017}
\begin{equation}\label{eq_16}
\arctan \left( x \right)=\,i\sum\limits_{m=1}^{\infty }{\frac{1}{2m-1}\left( \frac{1}{{{\left( 1+2i/x \right)}^{2m-1}}}-\frac{1}{{{\left( 1-2i/x \right)}^{2m-1}}} \right)}.
\end{equation}
Recently Jes\'us Guillera found a simple and elegant proof of this formula for the arctangent function (see \cite{Abrarov2017} for details). Despite simplicity the formula \eqref{eq_16} demonstrates a very rapid convergence especially when its argument $x$ tends to zero. Although the rational number ${{u}_{2}}$ requires a large number of digits in its numerator and denominator, nevertheless, due to relation
$$
\frac{1}{\left| {{u}_{2}} \right|}<<\frac{1}{{{u}_{1}}}
$$
its application provides more rapid convergence of the second arctangent function in the two-term Machin-like formula \eqref{eq_3} for pi. Consequently, the second term associated with rational number ${{u}_{2}}$ in the two-term Machin-like formula \eqref{eq_3} for pi requires a smaller truncating integer in computation of the constant pi. As one can see, despite large number of the digits in numerator and denominator of the value ${{u}_{2}}$, this computational approach may be advantageous in algorithmic implementation.

Substituting equation \eqref{eq_16} into the two-term Machin-like formula \eqref{eq_3} for pi we have
\[
\begin{aligned}
\frac{\pi }{4}=&\,i\sum\limits_{m=1}^{\infty }\frac{1}{2m-1}\left( {{2}^{k-1}}\left( \frac{1}{{{\left( 1+2i{{u}_{1}} \right)}^{2m-1}}}-\frac{1}{{{\left( 1-2i{{u}_{1}} \right)}^{2m-1}}} \right)\right.\\
&\left.+\frac{1}{{{\left( 1+2i{{u}_{2}} \right)}^{2m-1}}}-\frac{1}{{{\left( 1-2i{{u}_{2}} \right)}^{2m-1}}} \right)
\end{aligned}
\]
or
\begin{equation}\label{eq_17}
\begin{aligned}
\pi = & \,4i\sum\limits_{m=1}^{\infty}\frac{1}{2m-1}\left( {{2}^{k-1}}\left( \frac{1}{{{\left( 1+2i{{u}_{1}} \right)}^{2m-1}}}-\frac{1}{{{\left( 1-2i{{u}_{1}} \right)}^{2m-1}}} \right)\right.\\
&\left.+\frac{1}{{{\left( 1+2i{{u}_{2}} \right)}^{2m-1}}}-\frac{1}{{{\left( 1-2i{{u}_{2}} \right)}^{2m-1}}} \right),
\end{aligned}
\end{equation}

The computational test reveals that with $k=27$, ${{u}_{1}}=85445659$ and corresponding rational value ${{u}_{2}}$ (see equation \eqref{eq_15} above), the truncated series expansion \eqref{eq_17} provides $16$ digits of pi per term increment. This convergence rate if faster than that of the Chudnovsky formula for pi
$$
\frac{1}{\pi }=\frac{12}{\sqrt{{{640320}^{3}}}}\sum\limits_{k=0}^{\infty }{{{\left( -1 \right)}^{k}}\frac{\left( 6k \right)!}{{{\left( k! \right)}^{3}}\left( 3k \right)!}\frac{13591409+545140134k}{{{\left( {{640320}^{3}} \right)}^{k}}}}
$$
providing $15$ digits of pi per term increment in truncation \cite{Arndt2001}. Furthermore, in contrast to the Chudnovsky formula for pi the proposed formula \eqref{eq_17} consists of the rational numbers only. Any irrational number involved in computation requires all digits. Specifically, if the number pi is supposed to be computed up to one trillion digits, then any irrational number involved in computation must contain all trillion digits. Therefore, the absence of the irrational numbers may also be advantageous in the proposed formula \eqref{eq_17} for pi.

There are several iteration-based algorithms providing enhanced convergence in computing pi. For example, the Brent--Salamin algorithm doubles a number of the correct digits of pi at each iteration \cite{Borwein1987, Arndt2001} and, therefore, just $25$ iterations are sufficient to produce a value of pi correct to over $45$ million digits. More rapid iteration algorithm, discovered by Borwein brothers, quadruples the number of correct digits at each iteration \cite{Borwein1987, Arndt2001}. One of the most rapid algorithms reported by Borwein {\it{et al}.} \cite{Borwein1989} provides quintic convergence that multiplies a number of the correct digits of pi by factor of $6$ at each iteration step. However, despite tremendously rapid convergence these iteration-based algorithms require irrational numbers appearing over and over again at each consecutive step of iteration. Perhaps this is one of the main reasons explaining why the most recent records \cite{Yee} in computing digits of pi were achieved by using the Chudnovsky formula that needs only one irrational number.

Applying equation \eqref{eq_13} the two-term Machin-like formula \eqref{eq_3} for pi can also be expressed alternatively in trigonometric form of the argument $1/{{u}_{2}}$ of the arctangent function
\begin{equation}\label{eq_18}
\frac{\pi }{4}={{2}^{k-1}}\arctan \left( \frac{1}{{{u}_{1}}} \right)+\arctan \left( \frac{1-\sin \left( {{2}^{k-1}}\arctan \left( \frac{2{{u}_{1}}}{u_{1}^{2}-1} \right) \right)}{\cos \left( {{2}^{k-1}}\arctan \left( \frac{2{{u}_{1}}}{u_{1}^{2}-1} \right) \right)} \right)
\end{equation}
as a complete analog of the equation \eqref{eq_6}. Although iteration procedure described above is more efficient for computation, the equation \eqref{eq_18} may be convenient to verify the results by using user-friendly mathematical languages like Mathematica or Maple. 

The following is an example of the Mathematica code showing the convergence rate by using equations \eqref{eq_16} and \eqref{eq_18}:
\small
\begin{verbatim}
(* Define integer k *)
k = 27;

(* Define value u1 *)
u1 = 85445659;

(* Compute value u2 *)
u2 = (Cos[2^(k - 1)*ArcTan[(2*u1)/(u1^2 - 1)]])/(1 -
    Sin[2^(k - 1)*ArcTan[(2*u1)/(u1^2 - 1)]]);

(* Approximation for pi, M is the truncating integer *)
piApprox[M_] :=
  N[4*I*Sum[(1/(2*m - 1))*(2^(k - 1)*(1/(1 + 2*I*u1)^(2*m - 1)
      - 1/(1 - 2*I*u1)^(2*m - 1)) + 1/(1 + 2*I*u2)^(2*m - 1)
          - 1/(1 - 2*I*u2)^(2*m - 1)), {m, 1, M}], 10000] //Re

Print["Number of coinciding digits with pi"]
piDigits[M_] := Abs[MantissaExponent[Pi - piApprox[M]]][[2]]
M = 1;
While[M <= 20, Print["At M = ", M, 
    " the number of coinciding digits is ", piDigits[M]]; M++]
\end{verbatim}
\normalsize

It should be noted that the values in equation \eqref{eq_18}
$$
\sin \left( {{2}^{k-1}}\arctan \left( \frac{2{{u}_{1}}}{u_{1}^{2}-1} \right) \right)
$$
and
$$
\cos \left( {{2}^{k-1}}\arctan \left( \frac{2{{u}_{1}}}{u_{1}^{2}-1} \right) \right)
$$
are rationals when ${{u}_{1}}$ is a rational number. This follows from the equations \eqref{eq_13} and \eqref{eq_14}.
\\
\begin{figure}[ht]
\begin{center}
\includegraphics[width=23pc]{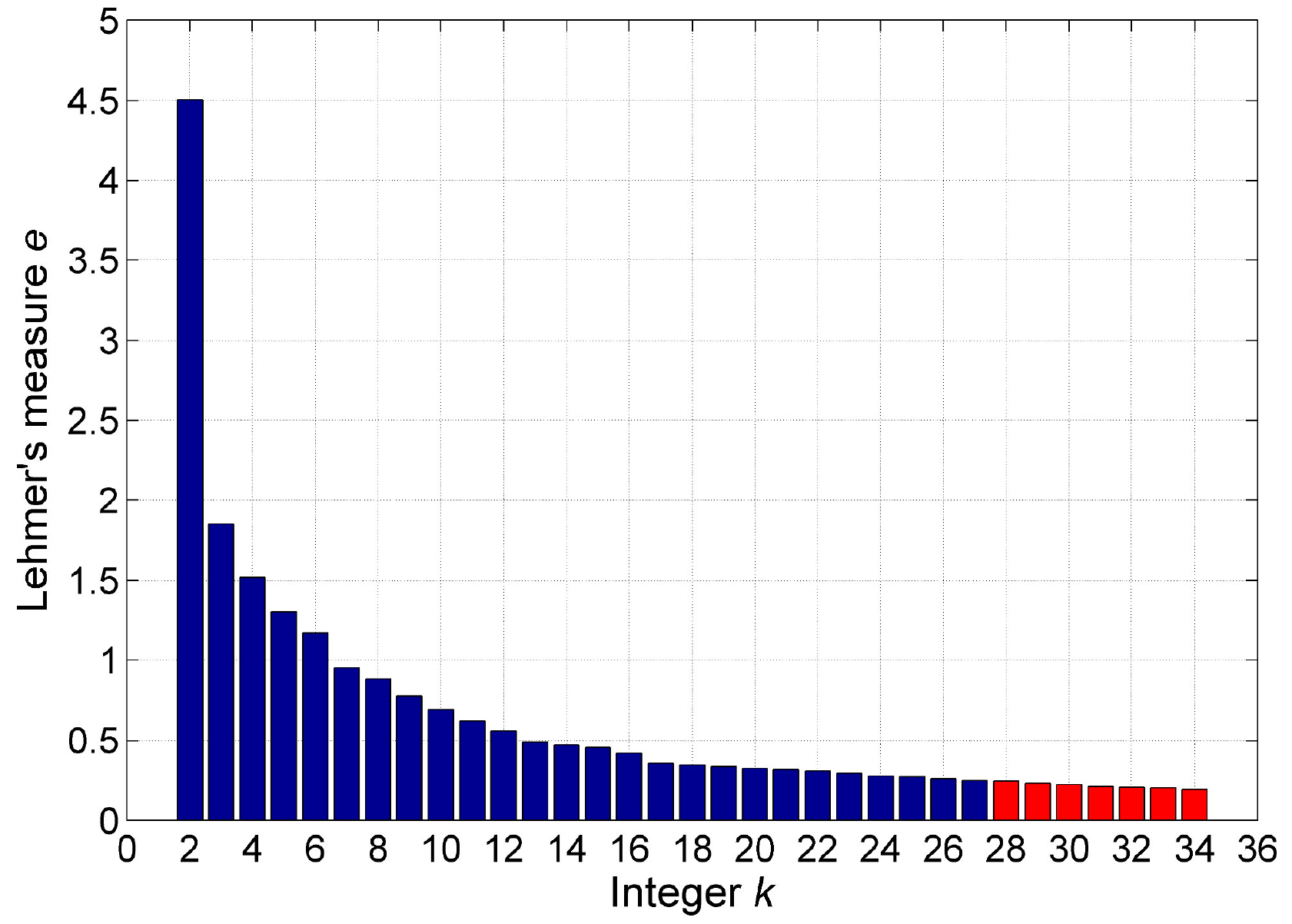}\hspace{1pc}%
\begin{minipage}[b]{28pc}
\vspace{0.2cm}
{\sffamily {\bf{Fig. 1.}} The Lehmer's measure $e$ as function of the integer $k$. The blue bars are computed with exact values of the rational integer $u_2$. The red bars are computed by using equation \eqref{eq_13}.}
\end{minipage}
\end{center}
\end{figure}

Figure 1 illustrates the Lehmer's measure $e$ dependence on the integer $k$ for the case when the rational value $u_1$ is computed according to the equation
$$
{{u}_{1}}=\left\lfloor \frac{{{a}_{m}}}{\sqrt{2-{{a}_{m-1}}}} \right\rfloor,
$$
where $\left\lfloor \ldots \right\rfloor$ denotes the floor function. The blue bars in Fig. 1 covering the range from $k = 2$ up to $k = 27$ are computed with exact values of the rational numbers $u_1$ while the red bars are computed by equation \eqref{eq_13}.

We found experimentally that the convergence rate is roughly equal to $4.1/e$. Consequently, by extrapolation we can expect that the convergence rate at $k=34$ can be increased up to $20$ digits per term increment once the rational value $u_2$ is found on a supercomputer. The application of the equation \eqref{eq_18} also predicts additional $20$ correct digits of pi per term increment by one.

\section{Computational methodology}

Although equation \eqref{eq_17} is simple and involves no surd numbers, it, nevertheless, can be simplified further in order to exclude all complex values in computation. This is possible to achieve by using the iteration technique described in this section.

Applying the following expansion coefficients
\[
\label{eq_19a}\tag{19a}
{{c}_{1}}\left( x \right)=1+2i/x
\]
and
\[
\label{eq_19b}\tag{19b}
{{c}_{m}}\left( x \right)={{c}_{m-1}}\left( x \right){{\left( 1+2i/x \right)}^{2}},
\]
one can express the series expansion for the arctangent function \eqref{eq_16} as
\setcounter{equation}{19}
\begin{equation}\label{eq_20}
\arctan \left( x \right)=i\sum\limits_{m=1}^{\infty }{\frac{1}{2m-1}\,\frac{{{{\bar{c}}}_{m}}\left( x \right)-{{c}_{m}}\left( x \right)}{{{c}_{m}}\left( x \right){{{\bar{c}}}_{m}}\left( x \right)}},
\end{equation}
where
$$
{{\bar{c}}_{m}}\left( x \right)=\operatorname{Re}\left[ {{c}_{m}}\left( x \right) \right]-i\operatorname{Im}\left[ {{c}_{m}}\left( x \right) \right]
$$
is the complex conjugate with respect to ${{c}_{m}}\left( x \right)$. Since
$$
{{\bar{c}}_{m}}\left( x \right)-{{c}_{m}}\left( x \right)=-2i\operatorname{Im}\left[ {{c}_{m}}\left( x \right) \right]
$$
and
$$
{{c}_{m}}\left( x \right){{\bar{c}}_{m}}\left( x \right)={{\operatorname{Re}}^{2}}\left[ {{c}_{m}}\left( x \right) \right]+{{\operatorname{Im}}^{2}}\left[ {{c}_{m}}\left( x \right) \right],
$$
the equation \eqref{eq_20} can be rearranged as
\begin{equation}\label{eq_21}
\arctan \left( x \right)=2\sum\limits_{m=1}^{\infty }{\frac{1}{2m-1}\,\frac{\operatorname{Im}\left[ {{c}_{m}}\left( x \right) \right]}{{{\operatorname{Re}}^{2}}\left[ {{c}_{m}}\left( x \right) \right]+{{\operatorname{Im}}^{2}}\left[ {{c}_{m}}\left( x \right) \right]}}
\end{equation}
Defining now
$$
{{a}_{m}}\left( x \right)=\operatorname{Im}\left[ {{c}_{m}}\left( x \right) \right]
$$
and
$$
{{b}_{m}}\left( x \right)=\operatorname{Im}\left[ {{c}_{m}}\left( x \right) \right],
$$
from equations \eqref{eq_19a} and \eqref{eq_19b} it follows that
$$
{{a}_{1}}\left( x \right)=2/x,
$$
$$
{{b}_{1}}\left( x \right)=1,
$$
$$
{{a}_{m}}\left( x \right)={{a}_{m-1}}\left( x \right)\left( 1-4/{{x}^{2}} \right)+4{{b}_{m-1}}\left( x \right)/x,
$$
$$
{{b}_{m}}\left( x \right)={{b}_{m-1}}\left( x \right)\left( 1-4/{{x}^{2}} \right)-4{{a}_{m-1}}\left( x \right)/x.
$$
Consequently, the equation \eqref{eq_21} can be conveniently rearranged as given by
\begin{equation}\label{eq_22}
\arctan \left( x \right)=2\sum\limits_{m=1}^{\infty }{\frac{1}{2m-1}\,\frac{{{a}_{m}}\left( x \right)}{a_{m}^{2}\left( x \right)+b_{m}^{2}\left( x \right)}}.
\end{equation}

The application of the Machin-like formulas \eqref{eq_1} for pi may represent a considerable interest in the context of present-day computational mathematics since the arctangent function can be expanded into series with very rapid convergence especially at smaller values of the Lehmer's measure. Moreover, since the integer $k$ in the equation \eqref{eq_3} can be in principle arbitrarily large, we have no any theoretical restrictions to reduce further the Lehmer's measure. For example, one of the rapid formulas is the well-known series expansion discovered by Euler \cite{Chien-Lih2005}
\begin{equation}\label{eq_23}
\arctan\left(x\right)=\sum\limits_{m=0}^{\infty}{\frac{{2^{2m}}{\left(m! \right)}^{2}}{\left(2m+1\right)!}\,\frac{x^{2m+1}}{{\left(1+x^2\right)}^{m+1}}}
\end{equation}
that provides very high-accuracy especially when the argument $x$ tends to zero. This tendency can be seen from the Fig. 2 showing how fast the error term \footnote{We imply the error term as a difference between the actual arctanget function and its truncated series expansion.} vanishes with decreasing the argument $x$ by absolute value with just $10$ summation terms in truncation. However, the numerical test reveals that the proposed series expansion \eqref{eq_22} of the arctangent function is more faster in convergence by many orders of the magnitude as it can be seen from the Fig. 3. Thus with only $10$ summation terms in truncation, at $x = 10^{-6}$ the equations \eqref{eq_23} and \eqref{eq_22} produce the error terms $2.7026 \times 10^{-127}$ and $4.54131 \times 10^{-134}$, respectively. Therefore, from these numerical examples we can see that it is more preferable to chose the series expansion \eqref{eq_22} of the arctangent function for computation of the constant pi.
\\
\begin{figure}[ht]
\begin{center}
\includegraphics[width=21pc]{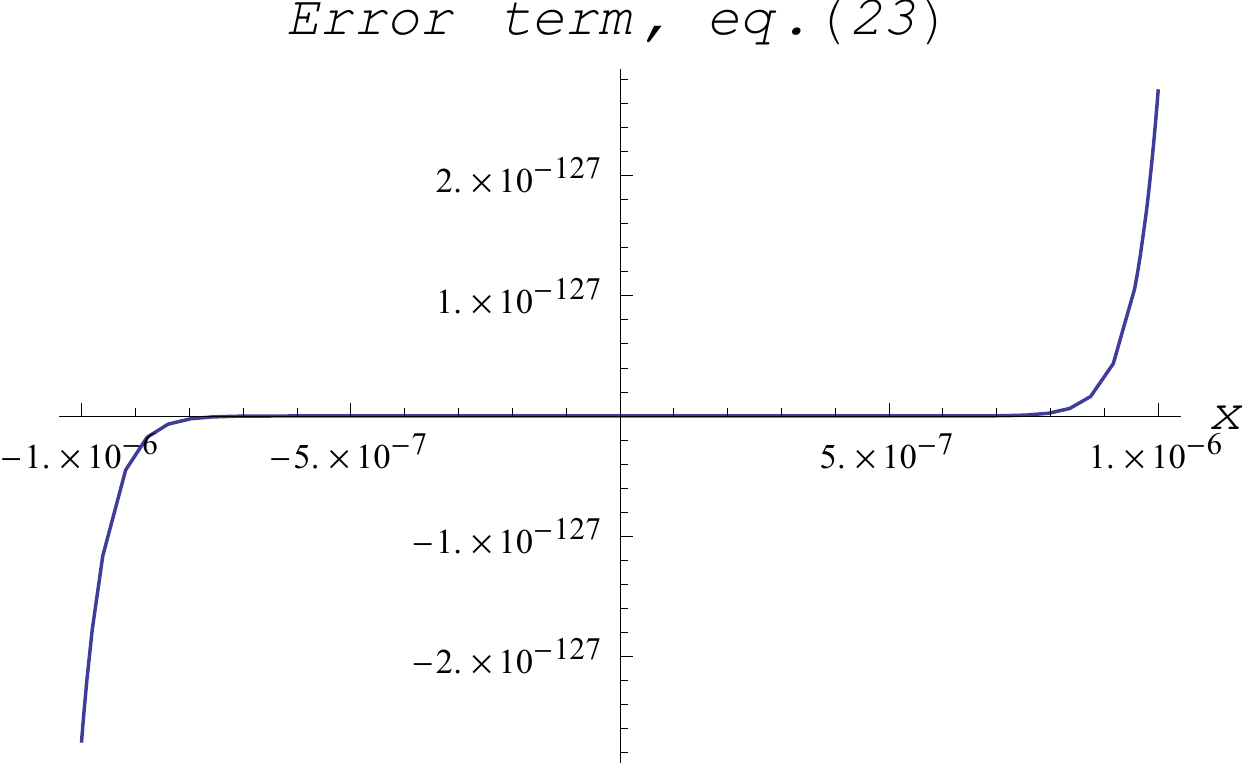}\hspace{1pc}%
\begin{minipage}[b]{28pc}
\vspace{0.2cm}
{\sffamily {\bf{Fig. 2.}} The error term inside the range $x \in \left[-10^{-6},10^6\right]$ at $10$ summation terms in truncation of the series expansion \eqref{eq_23}.}
\end{minipage}
\end{center}
\end{figure}
Substituting the equation \eqref{eq_22} into the two-term Machin-like formula \eqref{eq_3} for pi yields
\[
\frac{\pi }{4}=2\sum\limits_{m=1}^{\infty }{\frac{1}{2m-1}\left( {{2}^{k-1}}\frac{{{a}_{m}}\left( 1/{{u}_{1}} \right)}{a_{m}^{2}\left( 1/{{u}_{1}} \right)+b_{m}^{2}\left( 1/{{u}_{1}} \right)}+\frac{{{a}_{m}}\left( 1/{{u}_{2}} \right)}{a_{m}^{2}\left( 1/{{u}_{2}} \right)+b_{m}^{2}\left( 1/{{u}_{2}} \right)} \right)}
\]
or
\begin{equation}\label{eq_24}
\pi =8\sum\limits_{m=1}^{\infty }{\frac{1}{2m-1}\left( {{2}^{k-1}}\frac{{{\alpha }_{m}}}{\alpha _{m}^{2}+\beta _{m}^{2}}+\frac{{{\gamma }_{m}}}{\gamma _{m}^{2}+\theta _{m}^{2}} \right)},
\end{equation}
where the corresponding expansion coefficients can be found by iteration as follows
$$
{{\alpha }_{1}}=2{{u}_{1}},
$$
$$
{{\beta }_{1}}=1,
$$
$$
{{\alpha }_{m}}={{\alpha }_{m-1}}\left( 1-4u_{1}^{2} \right)+4{{\beta }_{m-1}}{{u}_{1}},
$$
$$
{{\beta }_{m}}={{\beta }_{m-1}}\left( 1-4u_{1}^{2} \right)-4{{\alpha }_{m-1}}{{u}_{1}},
$$
and
$$
{{\gamma }_{1}}=2{{u}_{2}},
$$
$$
{{\theta }_{1}}=1,
$$
$$
{{\gamma }_{m}}={{\gamma }_{m-1}}\left( 1-4u_{2}^{2} \right)+4{{\theta }_{m-1}}{{u}_{2}},
$$
$$
{{\theta }_{m}}={{\theta }_{m-1}}\left( 1-4u_{2}^{2} \right)-4{{\gamma }_{m-1}}{{u}_{1}}.
$$
As we can see, the equation \eqref{eq_24} is significantly simplified and excludes all complex numbers in computing pi. 
\\
\begin{figure}[ht]
\begin{center}
\includegraphics[width=21pc]{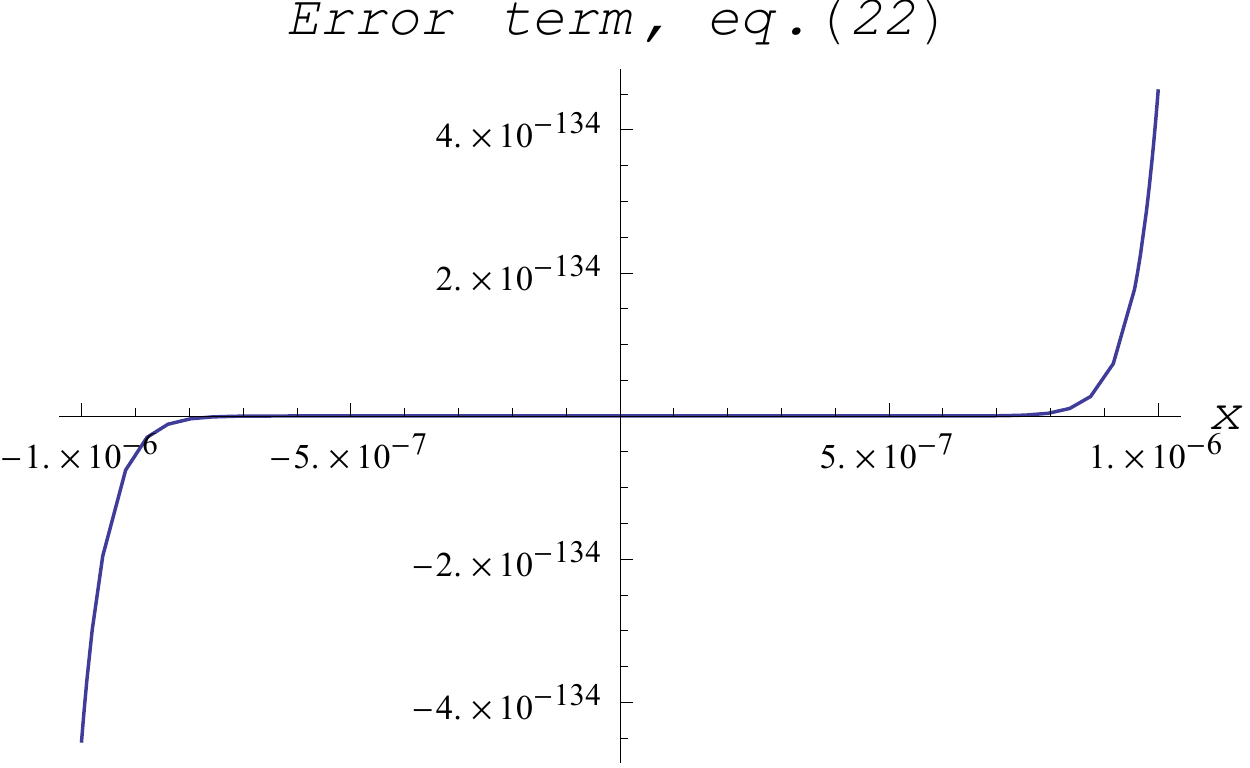}\hspace{1pc}%
\begin{minipage}[b]{28pc}
\vspace{0.2cm}
{\sffamily {\bf{Fig. 3.}} The error term inside the range $x \in \left[-10^{-6},10^6\right]$ at $10$ summation terms in truncation of the series expansion \eqref{eq_22}.}
\end{minipage}
\end{center}
\end{figure}

The following is an example of the Mathematica code implemented according to equation \eqref{eq_24}:
\small
\begin{verbatim}
(* Integer k *)
k = 6;

(* Rational number u1 *)
u1 = 40;

(* Rational number u2 *)
u2 = - 2634699316100146880926635665506082395762836079845121/
   38035138859000075702655846657186322249216830232319;

(* First set of the expansion coefficients *)
alpha[1] := alpha[1] = 2*u1
beta[1] := beta[1] = 1
alpha[m_] := alpha[m] = alpha[m - 1]*(1 - 4*u1^2) + 4*beta[m - 1]*u1
beta[m_] := beta[m] = beta[m - 1]*(1 - 4*u1^2) - 4*alpha[m - 1]*u1

(* Second set of the expansion coefficients *)
gamma[1] := gamma[1] = 2*u2
theta[1] := theta[1] = 1
gamma[m_] := gamma[m] = gamma[m - 1]*(1 - 4*u2^2) + 4*theta[m - 1]*u2
theta[m_] := theta[m] = theta[m - 1]*(1 - 4*u2^2) - 4*gamma[m - 1]*u2

(* Pi formula (24) *)
piApprox[M_] := 
  8*Sum[(1/(2*m - 1))*(2^(k - 1)*alpha[m]/(alpha[m]^2 + beta[m]^2) + 
       gamma[m]/(gamma[m]^2 + theta[m]^2)), {m, 1, M}];

(* Display pi with 100 decimal digits *)
Print["Actual value of pi is       ", N[Pi, 100]]
Print["Approximated value of pi is ", N[piApprox[25], 100]]
\end{verbatim}
\normalsize

\section*{Acknowledgment}
This work is supported by National Research Council Canada, Thoth Technology Inc. and York University.

\bigskip

\end{document}